\documentclass[11pt]{article}

\usepackage{amssymb}

\newcommand{\R}{\mathbb{R}}
\newcommand{\C}
     {\mathbb{C}}
\newcommand{\N}{\mathbb{N}}
\newcommand{\T}{\mathbb{T}}
\newcommand{\Z}{\mathbb{Z}}
\newcommand{\be}{\begin{equation}}
\newcommand{\ee}{\end{equation}}

\newtheorem{thm}{Theorem}[section]
\newtheorem{lem}[thm]{Lemma}

\newtheorem{prop}[thm]{Proposition}
\newtheorem{rmk}[thm]{Remark}

\begin{document}

\begin{center}
\noindent {\Large \bf{Global  attractor and asymptotic smoothing effects for the weakly damped cubic
Schr\"odinger equation in $L^2(\T)$}}
\end{center}
\vskip0.2cm
\begin{center}
\noindent
{\bf Luc Molinet}\\
{\small L.A.G.A., Institut Galil\'ee, Universit\'e Paris-Nord,\\
93430 Villetaneuse, France.} \vskip0.3cm
E-mail : molinet@math.univ-paris13.fr.
\end{center}
\vskip0.5cm \noindent {\bf Abstract.} {\small We prove that the weakly damped cubic  Schr\"odinger flow in
$L^2(\mathbb{T})$ provides a dynamical system that possesses a
global attractor. The proof relies on a sharp study of the behavior of the  associated flow-map
 with respect to the weak  $ L^2(\T) $-convergence inspired by \cite{L}.
 Combining the compactness in $ L^2(\T) $ of the attractor with the approach developed in \cite{G},
  we show that the attractor is actually a compact set of $ H^2(\T) $.
   This asymptotic smoothing effect is optimal in view of the  regularity of the steady states.
 }\vspace{2mm}\\

\section{Introduction}
The cubic nonlinear Schr\"odinger equation (NLS) can be derived as an asymptotic model to describe
 long wave propagation in different dispersive media. In some physical contexts, an exterior forcing and
    damping effects have to be taken into account and this can lead to the following
    cubic NLS equation
\begin{equation}\label{equation}
u_t+\gamma u +iu_{xx}\mp i|u|^2u=f,
\end{equation}
\noindent where $\gamma >0$ is the damping parameter and $ f$ is the forcing term.
 In this paper we focus on the case where $ u(t,x)$ is a function from $ \R_+\times \T $ to  $ \C $, with
  $ \T=\R/2\pi \Z $, and $ f\in L^2(\T) $ does not depend on time. Also since the sign in front of the 
   nonlinear term will not play any role in our analysis, we will take the $ + $ sign in all this paper.
   
    It is well-known since the work of Bourgain
   \cite{Bo1} that (\ref{equation}) provides an infinite dimensional dynamical system on $ H^s(\T) $
    for $ s\ge 0 $. Using an a priori estimate in $ H^1(\T) $ related to the energy conservation of the classical cubic NLS equation,
     the existence of a global attractor in $ H^1(\T) $ can be obtained by a standard method (see for instance \cite{Wang} or \cite{G1} where the additional regularity of the attractor is also proved). This method  principally contains 
      two steps. A first step consists in proving the continuity of the flow-map associated with the equation with respect to the
       weak topology of the phase space. This ensures the existence of a compact global attractor for the weak topology.
        The second step uses the argument of Ball \cite{Ball} to convert the weak convergence to the attractor into a strong one.
     The standard way to prove the first step is to use the well-posedness of the equation in a larger function space where
      the phase space is compactly embedded (cf. for instance \cite{Gh1}). This approach cannot be applied to (\ref{equation}) in
       $ L^2(\T) $ since the well-posedness of this equation is not known in  such a space. Actually, the
        strong ill-posedness of the classical cubic NLS equation in $H^s(\T) $, $s<0$,  has been even proved (cf. \cite{CCT2},
         \cite{L}).
          In \cite{OL} Goubet and the author used another approach involving the so called  Kato local smoothing
          effect for (\ref{equation}) on $ \R $ to establish the weak continuity of the flow-map in $ L^2(\R) $. 
           The situation  in $ L^2(\T) $ seems more complicated since, as was shown in \cite{L}, the flow-map of the classical cubic NLS equation
            is discontinuous for the weak topology of $ L^2(\T) $. Note however that the result in \cite{L} shows that the flow-map associated with the  modified Schr\"odinger equation (see (\ref{decomp} below) introduced by Christ \cite{C} is continuous for this weak topology. Let us mention here that the well-posedness of the  modified  NLS equation in a function space where $ L^2(\T) $ is continuously (but not compactly) embedded is obtained in \cite{C} (see also \cite{GH} for a related 
         result). Using  results of \cite{L}, 
                     we clarify
         some behaviors of the flow-map of (\ref{equation})with respect to the weak  $ L^2(\T) $-convergence.      
          From this  information supplemented with  the argument of Ball  we deduce the existence of a  global attractor in $L^2(\T) $. Finally, combining
           the approach developed in \cite{G1}-\cite{G} with the compactness of the attractor, we prove that the global attractor actually
            belongs to $H^2(\T) $ which can be viewed as  an asymptotic smoothing effect. This smoothing effect is optimal since for $ f$ belonging to $L^2(\T) $
             but  not to $ H^s(\T) $, with $ s>0$, the steady state to
            (\ref{equation}) does not belong to $ H^s(\T) $ for $ s>2$.

Denoting  by $ S(\cdot) $ the nonlinear  group associated with (\ref{equation}),  i.e.
 $
 S(t) u_0:= u(t) , \; t\in \R ,\;  $ where $ u $ is the solution of (\ref{equation}) associated with the initial data $ u_0$, 
 our main result is as follows :

\begin{thm}\label{attractor}
The nonlinear group $S(\cdot)$ associated with (\ref{equation}) provides an infinite-dimensional dynamical
system in $L^2(\T)$ that has a global attractor
$\mathcal{A}$ which is  a compact set of $ H^2(\T) $. More precisely, $\mathcal{A}$ is a connected and compact set of $ H^2(\T) $, invariant 
 (positively and negatively) by $ S(\cdot) $ that attracts for the  $ L^2(\T) $-metric  all positive  orbits uniformly with respect to bounded sets
  of initial data in   $ L^2(\T) $.

\end{thm}
\begin{rmk}
Exactly the same proof as in Section 5 below shows that the $ L^2(\R ) $ global attractor
  to (\ref{equation}) on the line, that was constructed in \cite{OL}, is  actually a  
 compact set of $ H^2(\R) $.
 \end{rmk}
The main new ingredient for proving Theorem \ref{attractor} is the following result on the behavior of the flow-map of  (\ref{equation}) with respect to the weak $ L^2(\T)$-convergence.
\begin{thm} \label{theoshro3}
Let $ \{u_{0,n}\} $ be a sequence of $ L^2(\T) $ converging
weakly to $ u_0 $ in $ L^2(\T) $ and let $ \{u_n\} $ be the
sequence of the associated emanating solutions
 of the weakly damped cubic Schr\"odinger equation (\ref{equation}).
  For any adherence value $ a_0 $ of $ \{ \| u_{0,n}\|_{L^2}^2\} $ there exists a continuous function $ t\mapsto a(t)$ from $ \R $ to $\R_+ $,
  with $ a(0)=a_0 $,
 and  a subsequence $ \{u_{n_k}\} $  of $\{u_n\} $ such that, for any $ t\in \R $,
 $ u_{n_k} (t) $ converges weakly in
 $ L^2(\T) $  to $ v(t) $ where $ v\in C(\R;L^2(\T))\cap L^4_{\mbox{loc}}(\R\times\T)$ is
 the unique solution  to
 \begin{equation}\label{modifiedequation}
\left\{
\begin{array}{l}
v_t +iv_{xx} +\gamma v + i |v|^2 v+  \displaystyle \frac{i}{\pi}\Bigl(a(\cdot)-
\| v(\cdot)\|_{L^2}^2\Bigr) v =f \\
v(0)=u_0
\end{array}
\right. \; .
\end{equation}
\end{thm}
\begin{rmk}
It is worth noticing that this theorem ensures that,  in sharp contrast with the case on the line (cf. \cite{OL}),  the flow-map associated with
 (\ref{equation}) is not continuous for the weak topology  of $ L^2(\T) $.
  Indeed, following \cite{L}, let $ u_0\in L^2(\T) $ be different from $ 0 $ and let $ \{\phi_n\}\subset L^2(\T) $
 be a sequence such that $ \phi_n \rightharpoonup 0 $ in $ L^2(\T) $
 and $ \| \phi_n\|_{L^2}^2 \to 2\pi $
   as $ n$ goes to infinity (one can take for instance $ \phi_n= e^{i n x} $). Setting $ u_{0,n}=u_0+
   \phi_n $, we get that $ u_{0,n} \rightharpoonup u_0 $ in $ L^2(\T) $ and
    $ \|u_{0,n} \|_{L^2}^2 \to \|u_0 \|_{L^2}^2+2\pi  $ as $
 n\to \infty $. On account of Theorem \ref{theoshro3}, the emanating solutions $ u_n $
  tend weakly in $ L^2(\T) $ for any fixed $ t\in\R$ to
 $ v $ satisfying (\ref{modifiedequation}). Observe  that $ w=v-u $ is solution of
 \begin{equation} \label{tr}
\left\{
\begin{array}{l}
w_t +iw_{xx} +\gamma w + i \Bigl(|v|^2 w +(\overline{w}u+w\overline{v}) u \Bigr)=-  \displaystyle \frac{i}{\pi}\Bigl(a(\cdot)-
\| v(\cdot)\|_{L^2}^2\Bigr) v \\
w(0)=0
\end{array}
\right. \; .
\end{equation}
Since $ v(0)=u_0 \neq 0 $ and $ a(0)=\|u_0 \|_{L^2}^2+2\pi \neq \| v(0)\|_{L^2}^2  $
 we infer that  the $ L^2(\T) $-norm of the right-hand side
 of (\ref{tr})  cannot vanish for small $ t\neq 0 $. Hence $w(t)=0 $ is not a solution of (\ref{tr}) and thus $ v(t)\neq u(t) $ for small
  $ t\neq 0 $.\\
 Finally, note that, since  $ L^2(\T) $ is compactly embedded in $ H^s(\T) $  for $ s<0 $, this proves that (\ref{equation}) is ill-posed in $ H^{s}(\T) $  as soon as $ s< 0 $.
  \vspace*{2mm}
 \end{rmk}
 This paper is organized as follows. In the next section we introduce some notation and the function spaces we will work with. Section 3 is
 devoted to the proof of Theorem \ref{theoshro3} and Section 4 is devoted to the existence of the global attractor. Finally in Section 5 we prove the asymptotic smoothing effect.
 \section{Function spaces and notation}
 When we affirm that a proposition is valid for $x+ $ (respectively $x-$) with $x \in \R $, we mean that there exists
  a small real number $ \epsilon>0 $ such that the proposition is valid for any real number in the interval
   $ ]x,x+\epsilon[ $ (respectively $]x-\epsilon,x[$). For $ (x,y)\in \R^2 $, $ x\lesssim y $
 means that there exists $C>0 $ such that $ x\le C y $. We will also denote by $ \varepsilon $  any function from $ \R_+ $ into itself
  that goes to zero at infinity.\\
 For  a $ 2\pi $-periodic function $ \varphi$, we define its space Fourier transform  by
$$
\hat{\varphi}(k):=\frac{1}{2\pi} \int_{\T } e^{-i k x} \,
\varphi(x) \, dx , \quad \forall k \in \Z \; ,
$$
and we denote by $ P_N \varphi $ and $  Q_N \varphi $ the $
L^2(\T) $  orthogonal projections on respectively the space  Fourier
modes $ |k|\le N $ and $ |k|>N $. \\
 We denote by $ V(\cdot) $ the
free group associated with the linearized Schr\"odinger equation,
$$
\widehat{V(t)\varphi}(k):=e^{-i k^2t} \,
\hat{\varphi}(k) , \quad k\in \Z \, .
$$
The Sobolev spaces $ H^s(\T) $ for $ 2\pi$-periodic functions are defined as usually and endowed with
$$
\|\varphi\|_{H^s(\T)}:=\|\langle k \rangle^{s}
\widehat{\varphi}(k) \|_{l^2(\Z)} =\|J^s_x \varphi
\|_{L^2(\T)} \, ,
$$
where $ \langle \cdot \rangle := (1+|\cdot|^2)^{1/2} $ and $
\widehat{J^s_x \varphi}(k):=\langle k \rangle^{s}
\widehat{\varphi}(k)
$. \\
For a function $ u(t,x) $ on $ \R\times\T $, we define its space-time Fourier transform by
$$
\hat{u}(\tau,\xi):={\mathcal
F}_{t,x}(u)(\tau,\xi):=\frac{1}{2\pi}\int_{\R\times\T} e^{-i (\tau
t+ k x)} \, u(t,x)  dt \, dx , \quad \forall (\tau,k) \in
\R\times\Z \, .
$$
 and define the  Bourgain spaces $ X^{b,s} $ and $ {\tilde X}^{b,s} $ of  functions on $\R\times \T $
  respectively endowed with the norm
$$
\| u \|_{X^{b,s} } : =
 \| \langle \tau+k^2\rangle^{b}  \langle k \rangle^s
  \hat{u}\|_{L^2(\R; l^2(\Z))}=
  \| \langle \tau\rangle^{b}  \langle k \rangle^s
  {\mathcal F}_{t,x}(V(-t) u ) \|_{L^2(\R; l^2(\Z))}\; .
$$
and
$$
\| u \|_{{\tilde X}^{b,s} } : =
 \| \langle \tau-k^2\rangle^{b}  \langle k \rangle^s
  \hat{u}\|_{L^2(\R; l^2(\Z))}=
  \| \langle \tau\rangle^{b}  \langle k \rangle^s
  {\mathcal F}_{t,x}(V(t) u ) \|_{L^2(\R; l^2(\Z))}\; .
$$
Finally,  for an open interval $ I\subset \R  $ we define the restriction in time   spaces $ X^{b,s}_I $ of functions
 on $ I \times \T $
 endowed with the norm
$$
\| u \|_{X_{I}^{b,s}}: =\inf_{v\in X^{b,s}} \{ \|
v\|_{X^{b,s}} , \, v(\cdot)\equiv u(\cdot) \hbox{ on } I
\, \}\;.
$$
It is worth noticing that the $ X^{b,s}_I $ spaces are Hilbert spaces with dual (for the $ L^2$-duality)  $ X^{-b,-s}_I $  and that for any $ \theta \in [0,1] $ it holds 
$$
X^{\theta b_1+(1-\theta)b_2,s }_I=[X^{b_1,s}_I, X^{b_2,s}_I]^\theta \, .
$$
Moreover, for $ b>1/2 $, $ X^{b,s}_I $ is
continuously embedding in $ L^\infty(I; H^s(\T)) $ with a constant
 of continuity that depends on $ b$ and on $ |I| $ the length of $ I $, i.e.
 \begin{equation} \label{Linfi}
\| u \|_{L^\infty(I; H^s(\T))} \le C(b,|I|) \| u \|_{X_{I}^{b,s}}, \; \forall u\in X^{b,s}_I
 \end{equation}
 \section{Proof of Theorem \ref{theoshro3}}
Theorem \ref{theoshro3} is based on the observation made in \cite{L} on the cubic NLS equation posed on the one-dimensional torus.
We first recall the following well-posedness result due to Bourgain (\cite{Bo1}) for
  (\ref{equation}). Let us mention that this result was established for the cubic Schr\"odinger equation
   without damping and forcing but the adaptations for  (\ref{equation}) are straightforward.
\begin{thm} \label{Boschro3}
Let $ s\ge 0$. For any $ u_0\in H^s(\T) $, $ f\in H^s(\T) $  and any $T>0 $, there exists a unique solution
$$
u\in   L^{4}(]-T,T[\times \T)
$$
satisfying  (\ref{equation}) in ${\mathcal D}'(]-T,T[\times \T)$. Moreover $ u \in C([-T,T];H^s(\T))\cap X^{1/2+,s}_{]-T,T[} $ and the map data to solution $ u_0\mapsto u $ is real analytic from $ H^s(\T) $ to $ C([-T,T];H^s(\T)) $.
\end{thm}
Let us recall that this theorem principally use the linear estimates in Bourgain's spaces for the free evolution and the retarded Duhamel operator
\begin{equation} \label{L1}
\| V(t) \varphi\|_{X^{b,s}_{]-T,T[}} \le C(T,b) \| \varphi \|_{H^s}\; , b\in \R ,\; s\in\R, \, 0<T<1, 
\end{equation}
and for any $ 0<\varepsilon <\!< 1 $ and $ 0<T<1 $,
\begin{equation} \label{L2}
\|\int_0^t V(t-t') g(t') \, dt' \|_{X^{b,s}_{]-T,T[}} \le C(b,\varepsilon)  \, T^{\varepsilon}
\|g\|_{X^{b-1+\varepsilon,s}_{]-T,T[}}\, , 1/2\le b< 1 ,
\end{equation}
 as well as the following linear  dispersive estimate 
 \begin{equation}
\|v\|_{L^4(\R\times \T)} \lesssim \|v\|_{X^{3/8,0}}, \quad \forall v\in X^{3/8,0}\, . \label{L4}
\end{equation}
This estimate is proven in 
 \cite{Bo1}
  for functions on $ \T^2 $ but also holds for functions on $ \R\times\T $ (See \cite{L1} for a shorter proof that  works also clearly on $ \R\times\T$ ). 
 Moreover, according to \cite{Gi}, (\ref{L4})
ensures that for $ 0<T<1 $ it holds
\begin{equation}
\|V(t) \varphi\|_{L^4(]-T,T[\times\T)} \lesssim T^{1/8} \|\varphi \|_{L^2(\T)},
 \quad \forall \varphi\in L^2(\T) \, , \label{L41}
\end{equation}
which gives directly the existence and uniqueness in $
L^4(]-T,T[\times\T)$ by classical $ TT^* $ arguments. On the
other hand, to prove that $ u \in X_{]-T,T[}^{1/2+,0} $ one has
to notice that, applying (\ref{L4}) with $ \overline{v} $,  (\ref{L4}) clearly also holds with $X^{3/8,0} $ replaced by $
{\tilde X}^{3/8,0} $. Therefore, 
  \begin{equation} \label{conju}
\|{\mathcal F}^{-1}( |\widehat{\overline{v}}|)\|_{L^4(\T^2)}
\lesssim \|\overline{v}\|_{{\tilde X}^{3/8,0}} =
\|v\|_{{X^{3/8,0}}}\;,  \quad  \forall v\in X^{3/8,0} .
  \end{equation}
 and   (\ref{L4})-(\ref{conju}) yield
\begin{equation} \label{L5}
\|u_1 u_2
\overline{u_3}\|_{X^{-1/2+,0}_{]-T,T[}}\lesssim
\prod_{i=1}^3 \|u_i\|_{X^{1/2,0}_{]-T,T[}} , \quad \forall u_i\in
X^{1/2,0}_{]-T,T[}\; .
\end{equation}
Writing the Duhamel formulation of (\ref{equation}), using  (\ref{L1})-(\ref{L2}) and  (\ref{L5}) and choosing
 some small positive real number $ \varepsilon $, one can eventually derive the key estimate:
 \begin{equation} \label{L6}
\|u\|_{X^{1/2+,0}_{]-T,T[}}\lesssim \|u_0\|_{L^2(\T)}
+T^{0+} \Bigl[
(\|u\|_{X^{1/2+,0}_{]-T,T[}}^2+\gamma)
\|u\|_{X^{1/2+,0}_{]-T,T[}}+\|f\|_{L^2(\T)} \Bigr] \; .
 \end{equation}
 This  leads to the local existence result in $X_{]-T,T[}^{1/2+,0} $.
Finally, the fact that the time of existence in Theorem \ref{Boschro3} can
be chosen arbitrarly large follows from the a priori bound on the
$ L^2(\T) $-norm of the solution (see (\ref{absorbingball})-(\ref{tnegatif}) below).
\vspace*{2mm}

Now, let $ u_0 \in  L^2(\T) $ and $ \{u_{0,n}\}\subset L^2(\T) $
be a sequence  converging weakly to $ u_0 $ in $ L^2(\T) $. Note
that,
 from Banach Steinhaus'theorem,  $
\{||u_{0,n}||_{L^2(\T)}\} $ is bounded in $ \R_+ $.
  It is well-known that  the solutions of (\ref{equation}), given by Theorem \ref{Boschro3},
satisfy for all $ t\in \R $,
\begin{equation}\label{energydiff}
\frac{1}{2}\frac{d}{dt} ||u||^2_{L^2(\T)}+\gamma ||u||^2_{L^2(\T)} =
{\rm Re} \; \int_{\T} f\bar{u}dx
\end{equation}
By Young's inequality and Gronwall's lemma, we deduce that
the $ L^2 $-solutions satisfy for  any $t \in \R_+ $,
\begin{equation}\label{absorbingball}
||u(t)||^2_{L^2(\T)}\le e^{-\gamma t}||u_0||^2_{L^2(\T)}+
\frac{1-e^{-\gamma t}}{\gamma^2}||f||^2_{L^2(\T)}.
\end{equation}
Performing the change of variables $ (t,x)\mapsto (-t,x) $ and
proceeding as above we also infer that for any  $t\in \R_-$ it
holds
\begin{equation}\label{tnegatif}
||u(t)||^2_{L^2(\T)}\le e^{3\gamma |t|}||u_0||^2_{L^2(\T)}+
\frac{e^{3\gamma |t|}-1}{\gamma^2}||f||^2_{L^2(\T)}
\end{equation}
Therefore, from  (\ref{energydiff}), we deduce
that for any $ (t_0,t_1)\in \R^2 $ with $ t_1>t_0 $,
\begin{eqnarray}\label{equicontinuity}
\Bigl| ||u(t_1)||^2_{L^2_x}-||u(t_0)||^2_{L^2_x}\Bigr| &= &
\Bigl|2\gamma \int_{t_0}^{t_1} ||u(\tau)||^2_{L^2(\T)}\, d\tau +
2{\rm Re} \int_{t_0}^{t_1}\int_{\T} \;  f\bar{u}(\tau) dx\,
d\tau  \Bigr| \nonumber \\
& \le  & |t_1-t_0|\Bigl[ 3 \gamma \Bigl( e^{3\gamma
|t_1|}||u_0||^2_{L^2(\T)}+ \frac{e^{3\gamma
|t_1|}-1}{\gamma^2}||f||^2_{L^2(\T)}\Bigr) \nonumber \\
 & & +\frac{1}{\gamma}||f||^2_{L^2(\T)} \Bigr]\; .
\end{eqnarray}
 Denoting by $ u_{n}$  the solution to
(\ref{equation}) associated with the initial data $ u_{0,n}$, this
last inequality ensures that the sequence $ \{t\mapsto
||u_{n}(t)||_{L^2(\T)}^2\} $ is uniformly equi-continuous on any
bounded interval of $ \R $. It follows from Ascoli's theorem that
there exists a subsequence $ \{t\mapsto
||u_{n_k}(t)||_{L^2(\T)}^2\} $ that converges to some function $
t\mapsto a(t) $ in $C([-T,T];\R_+) $ for any $ T>0 $. Moreover, from
Theorem \ref{Boschro3} we know
  that $ \{u_{n_k}\} $ is bounded in $X^{1/2+,0}_{]-T,T[} $
  and thus, up to the extraction of  a subsequence, converges weakly to some $ v$ in   $X^{1/2+,0}_{]-T,T[} $.\\
Now, in (\cite{L}, Lemmas 2 $\& $ 3 ) it is proven that the nonlinear term of the modified Schr\"odinger equation  introduced in \cite{C}: 
\begin{equation}\label{decomp}
\Lambda(u):=|u|^2u -\frac{1}{\pi} \|u\|^2_{L^2} u
\end{equation}
is continuous from $(X^{1/2+,0}_1)^3 $
into $ X^{-7/16,0}_1$ equipped with their respective weak topology.
 We thus rewrite  the Duhamel formulation for $ u_n
$ in the following way :
 \begin{eqnarray}
u_n(t) &= & V(t) u_{0,n} - i \int_0^t V(t-t') \Lambda(u_n(t')) \, dt'  \nonumber \\
 &  & -\frac{i}{\pi}\int_0^t  V(t-t')
 \Bigl( \| u_{n}(t')\|_{L^2}^2 u_n(t')\Bigr)\, dt' - \gamma \int_0^t  V(t-t')u_n(t') d t'\nonumber \\
 & & + \;\int_0^t  V(t-t')f d t'
   \; .
\end{eqnarray}
Since $ u_{n_k} \rightharpoonup v $ in $X^{1/2+,0}_{]-T,T[} \hookrightarrow C([-T,T]; L^2(\T)) $ and $ a_{n_k}(\cdot) \to a(\cdot) $ in 
 $C([-T,T]; L^2(\T)) $, it follows that $a_{n_k}(\cdot) u_{n_k} \rightharpoonup a(\cdot) v $ in  $C([-T,T]; L^2(\T)) $. 
According to the continuity of the Duhamel operator from  $C([-T,T]; L^2(\T)) $ into
 itself,  the linear estimates (\ref{L1})-(\ref{L2}), the continuity
 result on $ \Lambda $ for the weak topology   and the above convergence results, we can pass to the limit to obtain that
\begin{eqnarray*}
v(t) &= & V(t) u_{0} - i \int_0^t V(t-t') \Lambda(v(t')) \, dt'   \\
 &  & \hspace*{-12mm} - \frac{i}{\pi}\,   \int_0^t V(t-t') (a(t') v(t'))\, dt'
 - \gamma \int_0^t  V(t-t')v(t') d t'
 + \int_0^t  V(t-t')f d t' \;
\end{eqnarray*}
and $ v$ is solution of the following Cauchy problem on $ ]-T,T[
$:
\begin{equation} \label{PDE}
\left\{
\begin{array}{l}
v_t +v_{xx}+\gamma v  +i \Lambda(v)+ \displaystyle \frac{i}{\pi} a(\cdot) v =f
\; \mbox{ in } \; {\cal D}'(]-T,T[\times \T) \\
v(0)=u_0
\end{array}
\right. \; .
\end{equation}
Proceeding exactly as for  the cubic Schr\"odinger equation, it
is easy to prove that  this Cauchy
 problem is globally well-posed\footnote{Note that the $ L^2 $-norm is
  controlled on any bounded interval of $ \R $} in $ H^s(\T) $, $s\ge 0 $, with a solution belonging  for all $T>0 $ to
$$ C([-T,T];H^s(\T)) \cap L^4(]-T,T[\times \T) $$
 with uniqueness in $ L^4(]-T,T[\times \T)$. Therefore, there exists only one possible limit and thus
 the sequence $ \{u_{n_k}\} $, and not only a subsequence of it,  converges weakly  to $ v$ in $ X^{1/2+,0}_{]-1,1[} $.
Moreover, using the equation satisfied by the $ u_{n} $ and the
uniform bound in $   L^\infty(]-T,T[; L^2(\T))\cap
L^4(]-T,T[\times \T)$, it is easy to check that for any smooth $
2\pi $-periodic function $\phi $, the family $ \{t\mapsto
(u_{n_k}(t),\phi)_{L^2} \}
 $ is bounded in $ C([-1,1]) $ and uniformly equi-continuous on $ [-1,1] $. Ascoli's theorem then ensures that $ (u_{n_k}, \phi) $ converges to $(v,\phi)$ on $ [-1,1] $ and thus $ u_{n_k}(t) \rightharpoonup v(t)
  $ in $ L^2(\T) $ for all $ t\in [-1,1] $. By direct iteration  this clearly also holds for all $ t\in \R $.  \\
 Finally, according to (\ref{decomp}), $ v$ can be also characterized as the unique solution in $ L^4(]-T,T[\times \T) $
  to
\begin{equation}
\left\{
\begin{array}{l}
v_t +iv_{xx} +\gamma v +i |v|^2 v + \displaystyle \frac{i}{\pi}\Bigl(a(\cdot)- \| v(\cdot)\|_{L^2}^2\Bigr) v =f \\
v(0)=u_0
\end{array}
\right. \; .
\end{equation}
 \section{Existence of the global attractor}
 Let us denote by $ S(t) $ the nonlinear  group associated with (\ref{equation}),  i.e.
 $$
 S(t) u_0:= u(t) , \; t\in \R \; .
 $$
On account of Theorem \ref{theoshro3} and (\ref{absorbingball}), we infer that
 the ball of $L^2(\T) $,
 $$
X:= \Bigl\{v\in L^2(\T) ,\, \|v\|_{L^2(\T)}\le M_0:=2\frac{||f||_{L^2(\T)}}{\gamma}\, \Bigr\}
 $$
  is a global
absorbing set for the dynamical system under consideration 
 and that  $S(t)$ acts continuously
on $X$ . To prove that there exists a global attractor it 
suffices to check the relative compactness in $ L^2(\T) $ of
sequences of the type $\{S(t_n) b_n\} $ with  $ t_n\uparrow
+\infty $ and $\{b_n\}\subset X $.
  This is the aim of the following proposition.
\begin{prop} \label{orbitcompact}
For any sequences $ \{b_n\}\subset  X $ and $ \{t_n\}\uparrow +\infty $,
 the sequence $ \{S(t_n) b_n\} $
  has an adherence value in $L^2(\T) $. \end{prop}
{Proof .}
We  combine Theorem \ref{theoshro3}  with the famous J.
Ball's argument (see \cite{Ball}, \cite{Wang}, \cite{MRW}).
Let $ \{b_n\}\subset  X $ and  let $ \{t_n\} $ be a sequence of positive real numbers
 that goes to infinity. From (\ref{absorbingball}) the sequence
  $ \{S(t_n) b_n\}  $ remains bounded in $ L^2(\T) $  and thus, up to the extraction of a subsequence, converges weakly in $ L^2(\T) $ to some $ v_0$.
  According to Theorem \ref{theoshro3}
  there exists a subsequence    $ \{S(t_{n_k}) b_{n_k}\} $  and a continuous function
 $ t\mapsto a(t) $ from $ \R $ to $ \R_+ $  such that  the
   solutions emanating from $\{S(t_{n_k} b_{n_k} )\} $ converge  weakly in $ L^2(\T) $ for all $ t\in \R $ to $ v(t) $ where
   $ v$ is the unique  solution to
    \begin{equation}\label{modifiedequation2}
\left\{
\begin{array}{l}
v_t +iv_{xx} +\gamma v + i |v|^2 v+  \displaystyle \frac{i}{\pi}\Bigl(a(\cdot)-
\| v(\cdot)\|_{L^2}^2\Bigr) v =f \\
v(0)=v_0
\end{array}
\right. \; .
\end{equation}
From (\ref{energydiff}) we infer that for $\tau>0 $ fixed and $n_k $ large enough,
\be\label{23}
 \|S(t_{n_k}) b_{n_k}\|^2_{L^2(\T)}=e^{-2\gamma
\tau}\|S(t_{n_k}-\tau) b_{n_k}\|^2_{L^2(\T)} -2{\rm Re}\int_0^\tau
\int_{\T}e^{-2\gamma s}\overline{f}\, S(t_{n_k}-s) b_{n_k} dsdx
\ee

\noindent where $ \| S(t_{n_k}-\tau) b_{n_k}\|^2_{L^2(\T)} \le M_0^2 $ and, according to the weak convergence and the dominated convergence theorem,

\be\label{24}\lim_{n_k\rightarrow +\infty}2{\rm Re}\int_0^\tau
\int_{\T}e^{-2\gamma s}\overline{f} \, S(t_{n_k}-s) b_{n_k}
dxds=2{\rm Re}\int_0^\tau \int_{\T}e^{-2\gamma s}\overline{f}\,
v(-s) dxds \; . \ee

\noindent On the other hand, using  the energy identity for equation (\ref{modifiedequation2}) ,  we
 get
\be\label{25} ||v_0||^2_{L^2(\T)}=e^{-2\gamma
\tau}||v(-\tau)||^2_{L^2(\T)} -2{\rm Re}\int_0^\tau
\int_{\T}e^{-2\gamma s}\overline{f}\, v(-s)dx ds. \ee
 But since $ S(t_{n_k}-\tau)  b_{n_k}\rightharpoonup v(-\tau) $ in $ L^2(\T) $, it follows from
  (\ref{absorbingball}) that
  $$
  \| v(-\tau)\|_{L^2(\T)}^2\le M_0^2 \; .
  $$
\noindent Gathering the above three   equalities, we thus infer that for any fixed
 $ \tau>0 $, 
\be\label{26}\limsup_{n_k\to +\infty}||S(t_{n_k})  b_{n_k}||^2_{L^2(\T)}\leq
||v_0||^2_{L^2(\T)} +2e^{-2\gamma \tau}M_0^2,\ee which ensures
that $ S(t_{n_k})  b_{n_k} $ converges actually strongly to $ v_0$ in $
L^2(\T) $. This completes the proof of  Proposition \ref{orbitcompact}. \vspace*{2mm}\\
Proposition \ref{orbitcompact} ensures the existence of a compact global attractor in $ L^2(\T) $.
 More precisely, from classical arguments (see for instance the proof of Theorem 1.1 in \cite{Temam}),
 it follows  that the positively invariant connected closed set
 $$
 {\mathcal A}:=\omega(X)=\bigcap_{s>0}\overline{ \bigcup_{t>s} S(t)X}
$$
is non-empty and attracts any bounded set of $ L^2(\T) $. The compactness of ${\mathcal A} $ follows as well. Indeed, let
$\{a_n\} \subset {\mathcal A} $. Taking a sequence $ \{t_n\}\uparrow +\infty$ and setting $ b_n=S(-t_n) a_n $, we get that $
a_n=S(t_n) b_n $ with $ \{b_n\} \subset {\mathcal A}\subset X $ and thus $\{a_n\} $ has got an adherence point in $ L^2(\T) $. Finally, it is worth noticing that, by construction, 
 $  {\mathcal A} $ is also negatively invariant.

\section{Asymptotic smoothing effect}
In this section we prove that the global attractor lies actually  in $ H^2(\T) $ and is moreover compact in this space. Following  the approach developed
 in \cite{G},  we split the solution $u(t)=S(t)u_0 $ emanating from $ u_0 $ into two parts by
 setting\footnote{Recall that $P_N $ and $ Q_N $ are the
 projections on respectively the spatial Fourier modes $ |k|\le N
 $ and $ |k|> N $.}
  \begin{equation}\label{eqv}
  v_t+iv_{xx}+\gamma v +i|v|^2 v = f- i P_N (|u|^2 u)+ i P_N(|v|^2 v)
  \end{equation}
 \begin{equation}
  w_t+iw_{xx}+\gamma w  = 
    -iQ_N (|w|^2 w-2|w|^2 u-w^2\overline{u}) -iQ_N(2|u|^2w+u^2\overline{w}) \nonumber 
  \label{eqw2}
  \end{equation}
 with initial conditions
 \begin{equation}\label{condvw}
  v(0)=P_N (u_0) \;\mbox{ and } \; w(0)=Q_N(u_0)\; .
  \end{equation}
  \begin{rmk}
  Proceeding as for the equation (\ref{equation}) it is easy to check that, $ u\in X^{1/2+,0}_T$  and $ 
   f \in L^2(\T)$ being given, the Cauchy problems (\ref{eqv}) and (\ref{eqw2}) are locally well-posed in $ L^2(\T) $. Hence, there exists $ \alpha> 0 $ and  a unique solution $ v\in X^{1/2+,0}_{]-\alpha,\alpha[} $ of 
   (\ref{eqv}) and $ w\in X^{1/2+,0}_{]-\alpha,\alpha[} $ of (\ref{eqw2}). Actually we will see   in this section that $ w\in C(\R_+;L^2(\T)) $  and $ v\in C(\R_+;H^2(\T)) $. 
     \end{rmk}
In \cite{G}, Goubet introduced this decomposition for the weakly damped  KdV equation. A first step
 of his analysis consists in proving  that the high frequency part
 $ w(t) $ is decreasing to $ 0 $ in $ L^2(\T) $.
 This decay of $ \| w(t)\|_{L^2(\T)} $ , which is uniform for all $ u_0 $ in the absorbing ball,  is obtained by using  the dispersive damping effect on the high-high frequencies interactions that occurs for the  nonlinear
  part of the KdV equation above
  $ H^{-1/2}(\T) $. This is related to the fact that the associated Cauchy problem is well-posed in $ H^s(\T) $ for $ s\ge -1/2 $.
 For the cubic Schr\"odinger equation the situation is more delicate since as recalled in the introduction this equation is ill-posed
 in $ H^s(\T) $ for $ s<0 $. Actually, due to some resonant parts in the nonlinear term, there is no damping effect on  high-high-high  interactions. To overcome this difficulty we will work directly on the global attractor and
   use in a crucial way that we already proved that it is compact in $ L^2(\T) $. Note that the a priori compactness of the global attractor is not required in
    \cite{G} where the compactness of the attractor can be obtained as a consequence of the asymptotic behavior of $ v$ and $ w$.
    
    The second step of the analysis in \cite{G} consists in proving an uniform bound in $ H^3(\T) $ 
     on $ v $. This uniform estimate follows from an uniform  bound   in $ L^2(\T) $ on the time 
      derivative $ v_t $ of $ v$. To get this last  bound  the author uses  that, in view of the equation satisfied by $ v$,  the low 
       frequencies $ P_N v_t $ belongs to any $ H^s(\T) $, $ s\in \R $. We will not be able to use this 
        approach here  since for $ v\in L^2(\T) $, $ P_N (|v|^2 v) $ does not belong a priori to any $ H^s(\T) $. 
       Inspired by \cite{Tsu} we will instead introduce the auxiliary function $ z:=Q_N (v -g)$, where $ g $ is
        defined by $\widehat{g}(k):=\frac{\widehat{ f}(k)}{-ik^2+\gamma} $, and prove that $ z$    
         is uniformly bounded in $ H^2(\T) $.
          \vspace*{2mm} \\
    The key proposition to derive the regularity of the attractor is the following.
    \begin{prop}\label{prop2}
 There exist functions $ h $ and $ K   \, :\, \R_+\to \R_+ $ with $ \lim_{t\to+\infty}
  h(t)=0 $  such that for  all $ N>0 $ large enough and all $ u_0\in {\mathcal A} $ the function $ v$ and $w$ constructed in (\ref{eqv})-(\ref{condvw}) satisfy
  \begin{equation}\label{estwv}
  \|w(t)\|_{L^2(\T)} \le h(t) \; \mbox{ and } \|v(t)\|_{H^2(\T)}\le K(N) \;  , \, \forall t\in \R_+ \, .
  \end{equation}
    \end{prop}
    With Proposition \ref{prop2} at hand it is  straightforward to check that $ {\mathcal A} $ is embedded in $ H^2(\T) $.  Indeed, let $ a \in
  {\mathcal A} $  and $ \{t_n\} \uparrow+\infty $. For all $ n\in\N $ we can write $ a $ as
  $$
  a=S(t_n) S(-t_n)a =S(t_n) b_n
  $$
  with $ b_n =S(-t_n) a \in {\mathcal A}$. From Proposition \ref{prop2} it follows that, for any $ n\in \N $,  $ a$ 
   can be decomposed
   as $ a=v_n+w_n $ with $ \| v_n \|_{H^2(\T)}\le K $ and $ \| w_n\|_{L^2(\T)}\to 0 $ as $ n\to +\infty $. Therefore  $
   a\in H^2(\T) $ and $\| a \|_{H^2(\T)}\le K. $ Hence, there exists $ K> 0 $, such that  the following uniform bound
    holds on the attractor :
   \begin{equation} \label{boundH2}
    \| a \|_{H^2(\T)}\le K \, , \; \forall a\in {\mathcal A}.
    \end{equation}
 \subsection{Proof of Proposition \ref{prop2}}
 \subsubsection{Preliminaries}
The $ L^2(\T) $-compactness of ${\mathcal A}$ ensures the following uniform bound on the $ L^2(\T)$-norm of the high  frequency part
  to  the elements of ${\mathcal A}$.
  \begin{prop}\label{prop1}
There exists a function $ \varepsilon $ from $ \R_+ $ into itself  that goes to zero at infinity such that 
  \begin{equation} \label{estprop1}
  \| Q_{N} a \|_{L^2(\T)} \le \varepsilon(N)\, , \; \forall a\in {\mathcal A} \; .
  \end{equation}
    \end{prop}
    Thanks to this remark we will have to prove a damping effect only on terms of the form
     $ P_{N/2} u_1 P_{N/2}u_2 \overline{Q_{N}u_3} $. This is the aim of the following lemma :
     \begin{lem} \label{damping}
     Let $ I \subset \R $ be a bounded interval and
     let $ u_i\in X^{1/2,0}_I $, $ i=1,2,3 $. Then for $ \epsilon>0 $ small enough it holds
     \begin{equation}
     \| P_{N/2}u_1 P_{N/2}u_2 \overline{Q_{N}u_3} \|_{X^{-1/2+\epsilon,0}_I}
     \lesssim N^{-1/4+2\epsilon} \prod_{i=1}^3 \| u_i \|_{X^{1/2,0}_I} \; .
     \end{equation}
     \end{lem}
    {Proof. } We take extensions $ v_i $ of the $ u_i $'s such that $
    \| v_i \|_{X^{1/2,0}}\le 2 \| u_i \|_{X^{1/2,0}_I} $. By duality we have  to prove that
    $$
     \sup_{\| w\|_{X^{1/2-\epsilon,0}}=1} \Bigl| \Bigl(w,P_{N/2}v_1
      P_{N/2}v_2 \overline{Q_{N}v_3})\Bigr)_{L^2(\R\times \T)} \Bigr|
      \lesssim N^{-1/4+2\epsilon} \prod_{i=1}^3 \| v_i \|_{X^{1/2,0}}\; .
     $$
     It thus suffices to estimate
     \begin{eqnarray*}
     J =  \int_{\R^3} \sum_{(k_1,k_2,k_3)\in A(N)}
      |\widehat{w}(\tau,k)|
 |\widehat{v_1}(\tau_1,k_1)| |\widehat{v_2}(\tau_2,k_2)|
|\widehat{\overline{v_3}}(\tau_3,k_3)|\, d\tau_1\, d\tau_2\, d\tau_3
     \end{eqnarray*}
     where $ \tau=\tau_1+\tau_2+\tau_3 $, $ k=k_1+k_2+k_3 $ and
     $$
     A(N):=\{(k_1,k_2,k_3)\in \Z^3, \, |k_1|\le N/2, \, |k_2|\le N/2 \mbox{ and } |k_3|> N \,  \}\; .
     $$
   To do this we will use the famous resonance relation for the  Schr\"odinger equation.
      Setting $ \sigma=\tau+k^2, \, \sigma_1=\tau_1+k_1^2, \, \sigma_2=\tau_2+k_2^2 $
  and $ \tilde{\sigma}_3=\tau_3-k_3^2 $, it holds
  \begin{equation}
   \sigma-\sigma_1-\sigma_2-{\tilde \sigma}_3 =2(k_3+k_1)(k_3+k_2)  \label{resonant} \;.
   \end{equation}
  This ensures that on $ \R^3\times A(N) $,
  $$
  \max(|\sigma|,|\sigma_1|,|\sigma_2|,|{\tilde \sigma}_3|) \gtrsim N^2 \; . $$
   Therefore we get, thanks to (\ref{L4}) and (\ref{conju}),
\begin{eqnarray*}
J & \lesssim  &  N^{-1/4+2\epsilon} \int_{\R^3} \sum_{(k_1,k_2,k_3)\in A(N)}
      |\sigma|^{1/8-\epsilon}  |\widehat{w}(\tau,k)|
  |\sigma_1|^{1/8}  |\widehat{v_1}(\tau_1,k_1)|  \\
   & &  |\sigma_2|^{1/8} |\widehat{v_2}(\tau_2,k_2)|
  |\tilde{\sigma}_3|^{1/8} |\widehat{\overline{v_3}}(\tau_3,k_3)|\, d\tau_1\, d\tau_2\, d\tau_3 \\
& \lesssim &N^{-1/4+2\epsilon} \|\
{\mathcal F}^{-1}( |\sigma|^{1/8-\epsilon} |\widehat{w}|)\|_{L^4(\R\times\T)}
  \|{\mathcal F}^{-1}(  |\tilde{\sigma}|^{1/8}  |\widehat{\overline{v_3}}|)\|_{L^4(\R\times\T)}  \\
 & &  \prod_{i=1}^2 \|{\mathcal F}^{-1}( |\sigma|^{1/8} |\widehat{v_i}|)\|_{L^4(\R\times\T)} \\
& \lesssim & N^{-1/4+2\epsilon}\|w \|_{X^{1/2-\epsilon,0}} \prod_{i=1}^3 \|v_i \|_{X^{1/2,0}} \; .
\end{eqnarray*}
 This completes the proof of the lemma. \vspace{2mm} \\

We are now in position to prove the Proposition \ref{prop2}.
    Let $ u_0 \in {\mathcal A} $ we decompose $ u(t)=S(t)u_0 $ by 
    \begin{equation}
    u(t)=v(t)+w(t) \label{tutu}
   \end{equation}
    where $ v$ and $ w$ are defined as in (\ref{eqv})-(\ref{condvw}).
     Note that (\ref{eqw2}) can be rewritten as 
     $$
       w_t+iw_{xx}+\gamma w  = 
     -i Q_N (|u|^2 u)+ i Q_N(|v|^2 v)
     $$
     which clearly ensures that (\ref{tutu}) holds.
\subsubsection{Decay in time of $ w$}
 From  Theorem \ref{Boschro3}, (\ref{L6}) and the fact that $ u $ belongs to the attractor, we know that
   for all $t\in \R $ (recall that $ {\mathcal A} $ is positively and negatively invariant by the flow),
 \begin{equation} \label{estuniformu}
  \|  u\|_{X^{1/2+,0}_{]t-1 ,t+1 [}}\lesssim M_0\mbox{ and } \| Q_N u\|_{X^{1/2,0}_{]t-1,t+1[}}\lesssim
   \varepsilon(N) \, .
  \end{equation}
Since $ u\in X^{1/2+,0}_{]-T,T[} $ for any $ T>0 $, proceeding as
in Theorem \ref{Boschro3}  it is easy to prove that the Cauchy
problem
 for $w $ is locally well-posed in $ L^2(\T) $ and thus 
 $w\in C([-\alpha,\alpha]; L^2(\T)) $ for some $ \alpha>0 $.
Moreover,  proceeding as in the proof of Theorem \ref{Boschro3}, we get the following estimate on $ w$ for all $ t\in]-\alpha,\alpha[ $ and $ 0<\delta< \min(|t-\alpha|, |t+\alpha|) $,
\begin{equation} \label{est1w}
\| w\|_{X^{1/2+,0}_{]t-\delta ,t+\delta [}}\lesssim \|w(t)\|_{L^2(\T)} + \delta^{0+} \| w\|_{X^{1/2+,0}_{]t-\delta ,t+\delta[ }}\Bigl(\| w\|_{X^{1/2+,0}_{]t-\delta ,t+\delta [}}^2
+\| u\|_{X^{1/2+,0}_{]t-\delta ,t+\delta [}}^2+1\Bigr) \; .
\end{equation}
Assuming that $ \|w(t)\|_{L^2(\T)} $ is bounded by some constant $ A>0 $ on  $ [0,T]$ for some positive time $ T\in ]0,\alpha[ $,
 we deduce that there exists $ \delta_0=\delta_0(A)>0 $ such that
 for $ 0<\delta<\delta_0 $ small enough,
\begin{equation} \label{est2w}
\| w\|_{X^{1/2+,0}_{]t-\delta ,t+\delta [}}\lesssim \|w(t)\|_{L^2(\T)}, \quad \forall t\in [0,T] \;.
\end{equation}
From now on, we fix $ 0<\delta<\delta_0 $ such that (\ref{est2w}) holds.
From this last inequality and (\ref{Linfi}) we infer that
\begin{equation}\label{est3w}
\inf_{\tau\in ]t, t+\delta[} \| w(\tau)\|_{L^2(\T)}\gtrsim  \| w(t)\|_{L^2(\T)}, \; \forall t\in [0,T] \; .
\end{equation}
 Multiplying (\ref{eqw2}) with $ 2 \overline{w} $ and integrating over $ \T $ we get
\begin{eqnarray}
\frac{d}{dt} \| w\|^2_{L^2(\T) }+2\gamma  \| w\|^2_{L^2(\T) } &\le&
  2 \Bigl|  \int_{\T}   (2|w|^2 u-w^2\overline{u}) \overline{w}\Bigr|\nonumber \\
 & & +2\Bigl|  \int_{\T} u^2\overline{w}^2  \Bigr| \label{estw1} \; .
\end{eqnarray}
Integrating (\ref{estw1}) with respect to time we obtain  the
following estimate
 for any $ t\in [0,T] $,
\begin{eqnarray}
 \| w(t+\delta)\|_{L^2(\T) }^2 &\le & \| w(t)\|_{L^2(\T) }^2 e^{-\gamma \delta} - \gamma
 \int_t^{t+\delta} e^{-\gamma (t+\delta-s)}\| w(s)\|_{L^2(\T) }^2 \, ds  \nonumber \\
 &&+2  \Bigl|\int_t^{t+\delta} e^{-\gamma (t+\delta-s)}  \int_{\T}   (2|w|^2 u-w^2\overline{u}) \overline{w}\, ds \Bigr|\
 \nonumber \\
& & + 2 \Bigl|\int_t^{t+\delta} e^{-\gamma (t+\delta-s)}  \int_{\T}  u^2\overline{w}^2 \, ds\Bigr|\  \label{est4w}
\; .
\end{eqnarray}
From (\ref{est3w}) we infer that,
\begin{equation} \label{est5w}
 - \gamma
 \int_t^{t+\delta} e^{-\gamma (t+\delta-s)}\| w(s)\|_{L^2(\T) }^2 \, ds\le -C \, (1-e^{-\gamma\delta} )\| w(t)\|_{L^2(\T) }^2 \; .
\end{equation}
Let us estimate  now  the two last  time integrals in
(\ref{est4w}). To do this we will extensively use that, following \cite{GR}, for $h\in X^{-1/2+\alpha,0}_{]t,t+\delta[}$ and $g\in X^{1/2+\alpha,0}_{]t,t+\delta[} $
 with  $ 0<\alpha<\!<1$, it holds  
\begin{equation}
\Bigr|\int_{t}^{t+\delta} e^{-\gamma (t+\delta-s)} \int_{\T} h(s,x) \overline{g(s,x)} dx \, ds \Bigl|
\lesssim C(\delta,\alpha) \| h\|_{X^{-1/2+\alpha,0}_{]t,t+\delta[}} \| g\|_{X^{1/2+\alpha,0}_{]t,t+\delta[}}\label{tot}\; .
\end{equation}
Indeed, taking time extensions $ \tilde{h} $ and $ \tilde{g} $ of $ h $ and $g$ such that $ 
 \|\tilde{g}\|_{X^{-1/2+\alpha,0}}\le 2  \|g\|_{X^{-1/2+\alpha,0}_{]t,t+\delta[}} $ and $ \|\tilde{h}\|_{X^{1/2+\alpha,0}}\le 2  \| h \|_{X^{1/2+\alpha,0}_{]t,t+\delta[}} $,  
  we have
 \begin{eqnarray*}
 \Bigr|\int_{t}^{t+\delta} e^{-\gamma (t+\delta-s)} \int_{\T} h(s,x) \overline{g(s,x)} dx \, ds \Bigl|  &=  &
  \Bigr|\int_{\R} \int_{\T} \tilde{h}(s,x) \chi_{[t,t+\delta]} e^{-\gamma (t+\delta-s)} \overline{ \tilde{g}(s,x)} dx \, ds \Bigl| \\
  &\lesssim & \|\tilde{h}\|_{X^{-1/2+\alpha,0}}  \|\chi_{[t,t+\delta]} e^{-\gamma (t+\delta-s)} \tilde{g}\|_{X^{1/2-\alpha,0}}
  \end{eqnarray*}
  with 
  \begin{eqnarray*}
   \|\chi_{[t,t+\delta]} e^{-\gamma (t+\delta-s)} \tilde{g}\|_{X^{1/2-\alpha,0}} & \lesssim & \| \chi_{[t,t+\delta]} e^{-\gamma (t+\delta-s)}\|_{L^\infty}
   \|\tilde{g}\|_{X^{1/2-\alpha,0}}\\
   & & + \| \chi_{[t,t+\delta]} e^{-\gamma (t+\delta-s)}\|_{H^{1/2-\alpha} }  \|\tilde{g}\|_{L^\infty(\R; L^2(\T))} \\
   & \lesssim & C(\alpha,\delta)  \|\tilde{g}\|_{X^{1/2+\alpha,0}}
   \end{eqnarray*}
 With (\ref{tot}) at hand, we deduce from (\ref{L5}), (\ref{estuniformu}) and (\ref{est2w}) 
that
\begin{eqnarray} \label{uyuy}
I_1 & := &  \Bigl|  \int_t^{t+\delta} e^{-\gamma (t+\delta-s)}\int_{\T}   (2|w|^2 u-w^2\overline{u}) \overline{w}\, ds\Bigr| \nonumber \\
&  \lesssim  & \|  (2|w|^2 u-w^2\overline{u})  \|_{X^{-1/2+,0}_{]t,t+\delta[}} \|w\|_{X^{1/2+,0}_{]t,t+\delta[}}\nonumber\\
& \lesssim & \|w\|_{X^{1/2+,0}_{]t,t+\delta[}}^3\Bigl(  \|w\|_{X^{1/2+,0}_{]t,t+\delta[}}+  \|u\|_{X^{1/2+,0}_{]t,t+\delta[}} \Bigr)\nonumber\\
& \lesssim & \|w(t)\|_{L^2(\T)}^3\Bigl(  \|w(t)\|_{L^2(\T)}+  M_0\Bigr)
\end{eqnarray}
To estimate the last time integral we split it  into two parts in the following way:
\begin{eqnarray}
I_2& := &  \Bigl| \int_t^{t+\delta} e^{-\gamma (t+\delta-s)} \int_{\T}  u^2\overline{w}^2 \, ds\Bigr| \nonumber \\
& = &    \Bigl| \int_t^{t+\delta} e^{-\gamma (t+\delta-s)}  \int_{\T} (Q_{N/2} u(s))(Q_{N/2} u(s)+2P_{N/2} u(s))  (\overline{w}(s))^2 \, ds\Bigr| \nonumber\\
& & + \Bigl| \int_t^{t+\delta} e^{-\gamma (t+\delta-s)}\int_{\T} (P_{N/2} u(s))^2  (\overline{w}(s))^2\, ds\Bigr|\nonumber \\
&= & I_{21}+ I_{22} \; .\label{split}
\end{eqnarray}
To estimate $ I_{21}$  we proceed as above and use (\ref{estuniformu})  to get
\begin{eqnarray}
 I_{21} & \lesssim &  \|u\|_{X^{1/2+,0}_{]t,t+\delta[}}  \|Q_{N/2} u\|_{X^{1/2+,0}_{]t,t+\delta[}}  \|w\|_{X^{1/2+,0}_{]t,t+\delta[}}^2 \nonumber \\
 & \lesssim & M_0\,  \varepsilon(N/2)  \| w(t)\|_{L^2(\T) }^2  \; .
 \end{eqnarray}
 Finally, to estimate $ I_{22}$ we use Lemma \ref{damping} (recall that $ w=Q_N w$)  and (\ref{est2w}) to obtain
 \begin{eqnarray} \label{yuyu}
 I_{22} & \lesssim &  \| (P_{N/2} v(s))^2 \overline{w}\|_{X^{-1/2+,0}_{]t,t+\delta[}} \|w\|_{X^{1/2+,0}_{]t,t+\delta[}} \nonumber \\
 & \lesssim & \frac{M_0^2}{N^{1/4-}}  \|w\|_{X^{1/2+,0}_{]t,t+\delta[}}^2 \nonumber \\
 & \lesssim &   \frac{M_0^2}{ N^{1/4-}} \| w(t)\|_{L^2(\T) }^2 \; .
 \end{eqnarray}
 Gathering (\ref{est4w})-(\ref{yuyu}) we thus infer that for all $ t\in [0,T] $,
 \begin{eqnarray}\label{eqff}
  \| w(t+\delta)\|_{L^2(\T) }^2 & -  & \| w(t)\|_{L^2(\T) }^2 e^{-\gamma \delta} \nonumber \\
  & \le & \Bigl[
 C_1\,  \Bigl( \|w(t)\|_{L^2(\T)}(  \|w(t)\|_{L^2(\T)}+  M_0)+M_0\,  \varepsilon(N/2) +\frac{M_0^2}{ N^{1/4-}} \Bigl) \nonumber \\
  &  &-C_2 (1-e^{-\gamma \delta}) \Bigr] \| w(t)\|_{L^2(\T) }^2\; .
 \end{eqnarray}
 Since $ w(0)=Q_N(u_0) $, according to Proposition \ref{prop1}, we can choose $ N>0 $ large enough so that the right-hand side of the above inequality is negative at $ t=0 $. By direct
  iteration in time and (\ref{est3w}) we thus infer that
 \begin{equation} \label{fo}
    \|w(t)\|_{L^2(\T)}\lesssim e^{-\gamma t}  \|w(0)\|_{L^2(\T)} \lesssim e^{-\gamma t}
     \|Q_N u_0 \|_{L^2(\T)}\lesssim e^{-\gamma t}  \varepsilon(N) \, , \; \forall t\in [0,T] \;.
    \end{equation}
    In particular, $ \|w(t)\|_{L^2(\T)} $ is bounded by $ A=C \|w(0)\|_{L^2(\T)} $ on $ [0,T] $ and  from the  local well-posedness of (\ref{eqw2})  we infer that $ w\in 
    C(\R_+;L^2(\T)) $ and that 
     (\ref{fo}) holds actually for any $ T>0 $. This proves the first assertion of Proposition \ref{prop2}.
\subsubsection{Estimate on $Q_N v$}
First since $ u=v+w $ we deduce from the preceding subsection that, for $ N $ large enough, $ v$ is well defined for all positive time and $ v\in C(\R_+;L^2(\T))$. 
Now, since by construction $ P_N v=P_N u $ and $ u $ belongs to the global attractor, we get thanks
 to (\ref{estuniformu}) that
\begin{equation}\label{eq1vZ}
\| P_N v(t) \|_{H^2} \lesssim \|P_N u(t)\|_{H^2} \le C(N), \, \forall t\ge 0 \, .
\end{equation}
It thus remains to control the high frequencies of $ v$. Inspired by \cite{Tsu} we introduce the functions
 $ g$  and $  g_N $ defined by
 \begin{equation}\label{defg}
\hat{g}(k):=\frac{\widehat{ f}(k)}{-ik^2+\gamma} \mbox{ and } g_N:= Q_N  g 
\end{equation}
so that $ g_N$  satisfies the equation
$$
\partial_t g_N +i \partial_{xx} g_{N} +\gamma g_N=Q_N f \; .
$$
Therefore, setting $ z:=Q_N v -g_N $, $z=Q_N z $ and is solution of
\begin{equation} \label{eqz}
\left\{
\begin{array}{l}
z_t +iz_{xx} +\gamma z +i Q_N(|v|^2 v)=0 \\
z(0)=-g_N
\end{array}
\right. \; .
\end{equation}
We plan to prove that $ z(t) $ is uniformly bounded in $ H^2(\T) $ for positive times. 
 We will need the following result on the behavior
 of $ g_N $ with respect to $ N $.
 \begin{lem}\label{lemestg}
 $g_N \in H^2(\T) \cap X^{1/2,1}_{]-1,1[}  $ and it holds
 \begin{equation}\label{estg}
 \|g_N\|_{H^2(\T)}+\|g_N\|_{X^{1/2,1}_{]-1,1[} }  \le \varepsilon(N)
  \end{equation}
where $ \varepsilon(N) \to 0 $ as $ N\to +\infty $.
  \end{lem}
  {\it Proof. }
  It is clear that
$$
\|g_N\|_{H^2(\T)}\le \|Q_N f \|_{L^2(\T)} \to 0 \mbox{ as } N\to +\infty  \; .
$$
Let now $ \psi\in C^\infty_0(]-2,2[)$ such that $ \psi\equiv 1 $ on $ [-1,1] $. It holds
\begin{eqnarray*}
\|g_N\|_{X^{1/2,1}_{]-1,1[}}  & \le  & \|\psi g_N \|_{X^{1/2,1}}= \|\langle \tau+k^2\rangle^{1/2}\langle k \rangle
\hat{\psi} \hat{g_N} \|_{L^2(\R\times \Z)} \\
& \le & \|\langle \tau \rangle^{1/2} \hat{\psi} \|_{L^2(\R)}
\|\langle k \rangle\hat{g_N} \|_{L^2(\Z)}+ \|\psi \|_{L^2(\R)}
\|\langle k \rangle^{2} \hat{g_N} \|_{L^2(\Z)} \lesssim \|g_N\|_{H^2(\T)} \; .\end{eqnarray*}
This completes the proof of the lemma. \vspace*{2mm}\\
It is worth noticing that combining (\ref{estprop1}),  (\ref{estg}),
 (\ref{estuniformu}), (\ref{est2w})  and (\ref{fo}), there exists $ \delta_0>0 $ such that 
\begin{equation}\label{estz0}
\|z(t)\|_{L^2(\T)} \le \varepsilon(N) \mbox{ and }
\|z\|_{X^{1/2+,0}_{]t- \delta_0,t+ \delta_0[} }  \lesssim M_0 ,\;  \forall t\ge 0,
\end{equation}
where  $ \varepsilon(N) \to 0 $ as $ N\to +\infty $. Therefore, taking $\beta>0 $ small enough, it holds
\begin{equation} \label{estuniformz}
\|z\|_{X^{1/2,0}_{]t- \delta_0,t+ \delta_0[} } \lesssim  \|z\|_{X^{0,0}_{]t- \delta_0,t+ \delta_0[} }^\frac{2\beta}{1+2\beta}
 \|z\|_{X^{1/2+\beta,0}_{]t- \delta_0,t+ \delta_0[} }^\frac{1}{1+2\beta} \lesssim \varepsilon'(N)
\end{equation}
where  $ \varepsilon'(N) \to 0 $ as $ N\to +\infty $.

According to the linear estimates (\ref{L1})-(\ref{L2}), to prove that the equation (\ref{eqz}) is globally well-posed in $ H^2(\T) $, it suffices to prove the following estimate :
\begin{lem} \label{lemestit}
Assuming that $ z\in X^{1/2+,2}_I $ for some time interval $ I\subset \R $ with $ |I|\le 1$. 
The following estimate holds :
\begin{equation} \label{estit}
\Bigl\| Q_N(|v|^2 v)\Bigr\|_{X^{-1/2+\epsilon,2}_I} \lesssim C(N)+\|z\|_{X^{1/2+,2}_I}
\end{equation}
\end{lem}
{\it Proof. } We decompose $ v$ as $ v=P_N u +z +g_N $ so that we have to estimate
$$
\Bigl\| Q_N\Bigl(|P_N u+z+g_N|^2 (P_N u+z+g_N)\Bigr)\Bigr\|_{X^{-1/2+\epsilon,2}_I}\: .
$$
Let us first estimate the expression  containing $ g_N$, i.e. terms of the form
$
\|Q_N(\overline{g_N} w_1 w_2) \|_{X^{-1/2+\epsilon,2}_I} $ or $\|Q_N(\overline{w}_1 g_N w_2) \|_{X^{-1/2+\epsilon,2}_I} $
 with $(w_1,w_2)\in\{g_N,P_N u, z\}^2 $. By the triangle inequality we can write
 \begin{eqnarray}
\|Q_N(\overline{g_N} w_1 w_2) \|_{X^{-1/2+\epsilon,2}_I} & \le & \|Q_N(\overline{g_N} w_1 w_2) \|_{X^{-1/2+\epsilon,0}_I}
+\|Q_N(\overline{g_N}  D_x^2\tilde{ w}_1 w_2 )\|_{X^{-1/2+\epsilon,0}_I} \nonumber\\
 & & +
\|Q_N(\overline{g_N} w_1 D_x^2 \tilde{w}_2) \|_{X^{-1/2+\epsilon,0}_I}
+\|Q_N(\overline{D_x^2 g_N} w_1 w_2) \|_{X^{-1/2+\epsilon,0}_I}  \nonumber \\
& & +\|Q_N(\overline{g_N}D_x^2 g_N  w_2) \|_{X^{-1/2+\epsilon,0}_I} \label{etoile}
\end{eqnarray}
with $(\tilde{w_1}, \tilde{w}_2)\in \{P_N u, z\} $. 
The  terms containing no derivative on $ g_N$ of the above right-hand side can be estimated thanks to (\ref{L4}),
 (\ref{estg}) and (\ref{estz0})  by
\begin{eqnarray}
& &\|g_N\|_{X^{1/2,0}_I} (\|P_N u\|_{X^{1/2,2}_I} +\|z\|_{X^{1/2,2}_I}+\|g_N\|_{X^{1/2,0}_I})
(\|P_N u\|_{X^{1/2,0}_I}+ \|z\|_{X^{1/2,0}_I}+\|g_N\|_{X^{1/2,0}_I})
 \nonumber \\
& &\hspace*{8mm} \lesssim\varepsilon(N)  (M_0^2 N^2 +M_0 \|z\|_{X^{1/2,2}_I}+\varepsilon(N))
 \; . \label{estit2}
\end{eqnarray}
For  the terms that contains two  derivatives on $ g_N $ we write
\begin{eqnarray}
\|Q_N(D_x^2 \overline{g_N} w_1 w_2) \|_{X^{-1/2+\epsilon,0}_I} & \lesssim &
 \|g_N\|_{H^2(\T)} \|w_1\|_{L^\infty_I L^\infty(\T)} \|w_2\|_{L^\infty_I L^\infty(\T)}\nonumber \\
 & \lesssim & \varepsilon(N)\|w_1\|_{L^\infty_I L^2(\T)}^{3/4} \|w_2\|_{L^\infty_t L^2(\T)}^{3/4}
 \|w_1\|_{L^\infty_I H^2(\T)}^{1/4} \|w_2\|_{L^\infty_I H^2(\T)}^{1/4}\nonumber \\
  & \lesssim & \varepsilon(N) (M_0^2 N^2+\varepsilon(N)\|z\|_{X^{1/2+,2}_I} +\varepsilon(N))
  \label{estit3} \; .
\end{eqnarray}
The terms of the form $\|Q_N(\overline{w} g_N w_2) \|_{X^{-1/2+\epsilon,2}_I} $ can be treated exactly 
 in the same way.
It remains to consider the terms where $ g_N$ is not involved.  From (\ref{L4}) and (\ref{estuniformz}) ,
\begin{equation} \label{estit4} 
\||z|^2 z \|_{X^{-1/2+\epsilon,2}_I} \lesssim \|z\|_{X^{1/2,0}_I}^2\|z\|_{X^{1/2,2}_I}
 \lesssim \varepsilon(N)^2 \|z\|_{X^{1/2,2}_I} ,
\end{equation}
and
\arraycolsep1pt
\begin{eqnarray}
\|P_N u  |z|^2 \|_{X^{-1/2+\epsilon,2}_I} & +  & \|\overline{P_N  u} z^2 \|_{X^{-1/2+\epsilon,2}_I}\nonumber \\
 & \lesssim & \|P_N u\|_{X^{1/2,0}_I} \|z\|_{X^{1/2,0}_I} \|z\|_{X^{1/2,2}_I}
+ \|P_N u\|_{X^{1/2,2}} \|z\|_{X^{1/2,0}}^2 \nonumber\\
 & \lesssim & \varepsilon(N) M_0 \|z\|_{X^{1/2,0}_I}+N^2 M_0 \varepsilon(N) N^{-2}\|z\|_{X^{1/2,2}_I} \nonumber\\
 & \lesssim & \varepsilon(N) M_0 \|z\|_{X^{1/2,2}_I}  \label{estit5} \; .
\end{eqnarray}
To deal with $ (P_N u)^2  \overline{z}$ we decompose it as 
 \begin{equation} \label{estit6}
 (P_N u)^2 \overline{z}   =  
   (Q_{N/2} P_N u)( Q_{N/2} P_N u+2 P_{N/2} u) \overline{z} 
  +  ( P_{N/2} u)^2  \overline{z}  :=A_1+A_2 \; .
  \end{equation}
  Clearly, (\ref{estuniformu}) yields
   \begin{equation} \label{estit7}
 \|  A_1 \|_{X^{-1/2+\epsilon,2}_I} \lesssim \varepsilon(N/2) M_0 \|z\|_{X^{1/2,2}_I}
  \; 
  \end{equation}
and using Lemma \ref{damping} it is easy to check that 
  \begin{equation} \label{estit71}
 \|  A_2 \|_{X^{-1/2+\epsilon,2}} \lesssim N^{-1/4}  M_0^2  \|z\|_{X^{1/2,2}_I} \; . 
  \end{equation}
Finally,
\begin{eqnarray}
\||P_N u|^2  z| \|_{X^{-1/2+\epsilon,2}} 
& \lesssim & \|P_N u\|_{X^{1/2,0}_I} \|P_N u \|_{X^{1/2,2}_I} \|z\|_{X^{1/2,0}_I}
+ \|P_N u\|_{X^{1/2,0}_I}^2 \|z\|_{X^{1/2,2}_I}\nonumber \\
& \lesssim  &  M_0^2 (N^2 N^{-2}+1)\|z\|_{X^{1/2,2}_I}\lesssim M_0^2 \|z\|_{X^{1/2,2}_I} \label{estit9} 
\end{eqnarray}
and
\begin{equation} \label{estit8}
\||P_N u|^2  P_N u | \|_{X^{-1/2+\epsilon,2}_I} \lesssim N^2 M_0^3 \; .
\end{equation}
Gathering all the above estimates, (\ref{estit}) follows. \vspace*{2mm} \\
From the above lemma  and (\ref{L1})-(\ref{L2}) we
 deduce that $ z \in X^{1/2+,2}_{]t-\delta,t+\delta[} $    for any $ t\ge
0 $ and any $ 0<\delta<  \delta_0 $. Moreover,  it holds
 $$
 \| z\|_{X^{1/2+,2}_{]t-\delta ,t+\delta[ }}\lesssim \|z(t)\|_{H^2(\T)} + C \, \delta ^{0+} \Bigl( C(N)+\|z \|_{X^{1/2+,2}_{]t-\delta ,t+\delta[ }}\Bigr) \; .
  $$
 This ensures that for $ \delta_0>0 $ small enough,
 \begin{equation} \label{estZ'}
\| z\|_{X^{1/2+,2}_{]t-\delta_0 ,t+\delta_0 [}}\lesssim \|z(t)\|_{H^2(\T)} +  C(N) \; .
\end{equation}
 We will proceed  as in the preceding
subsection. From now on we fix $  0<\delta < \delta_0 $
such that (\ref{estZ'}) holds. As in (\ref{est3w}), this implies
that
\begin{equation}
\inf_{\tau\in ]t-\delta, t+\delta[} \| z(\tau)\|_{H^2(\T)}\ge C \| z(t)\|_{H^2(\T)}-C(N), \; \forall t\ge 0  \; .
\end{equation}
On the other hand,
taking the real part of the $ H^2(\T) $ hermitian-product of (\ref{eqz}) with $ 2 z$, we get
\begin{eqnarray}
\frac{d}{dt} \| z\|_{H^2(\T) }^2+2\gamma  \| z\|_{H^2(\T) }^2 & = & \Re \Bigl[-i 2(\!(Q_N(|v|^2 v) , z)\!)_{H^2}
\Bigr] \;  .
\end{eqnarray}
Integrating with respect to time this implies the following estimate for any $ t\ge 0 $,
\begin{eqnarray} \label{eto}
 \| z(t+\delta)\|_{H^2(\T) }^2 &\le & \| z(t)\|_{H^2(\T) }^2 e^{-\gamma \delta} - \gamma
 \int_t^{t+\delta} e^{-\gamma (t+\delta-s)}\| z(s)\|_{H^2(\T) }^2 \, ds  \nonumber \\
 &&+ 2  \Bigl[ \int_t^{t+\delta} e^{-\gamma (t+\delta-s)} \Im \Bigl[(\!(Q_N(|v(s)|^2 v(s)) , z(s))\!)_{H^2}
  \Bigr]\, ds \nonumber \\
   &\le & \| z(t)\|_{H^2(\T) }^2 e^{-\gamma \delta} +(1-e^{-\gamma \delta}) \Bigl(C(N)-
  C  \| z(t)\|_{H^2(\T) }^2 \Bigr)\nonumber \\
 &&+ 2\Bigl|   \Im \Bigl[  \int_t^{t+\delta} e^{-\gamma (t+\delta-s)}(\!(|v(s)|^2 v(s) , z(s))\!)_{H^2}\, ds
  \Bigr]\Bigr|
\end{eqnarray}
To estimate the last term of the above right-hand side we decompose $ v$ as in Lemma 
\ref{lemestit}. In view of  (\ref{tot}) and (\ref{etoile})-(\ref{estit9}) to get the following estimate :
\begin{equation}\label{est54}
\Bigl|   \Im \Bigl[  \int_t^{t+\delta}e^{-\gamma (t+\delta-s)} (\!(|v(s)|^2 v(s) , z(s))\!)_{H^2}\, ds
  \Bigr]\Bigr| \lesssim 
 C(N)+ \varepsilon(N) \| z\|_{X^{1/2+,2}}^2 \; ,
\end{equation}
we only have to care about 
$$
I:=\Bigl| \int_t^{t+\delta} e^{-\gamma (t+\delta-s)} \Im  (\!(|P_N u(s)|^2 z(s),z(s)\!  )_{H^2}\, ds\Bigr| \; .
$$
To deal with this term we decompose $\Im  (\!(|P_N u(s)|^2 z(s),z(s)  )\!)_{H^2} $ as 
$$
 \Im (\!(|P_N u|^2 z,z  )\!)_{H^1}\, ds\ 
 + \Im \int_{\T}\Bigl(2\partial_x( |P_N u|^2)   \partial_x z+
 \partial_x^2( |P_N u|^2)     z \Bigr)   \partial_x^2 \overline{z}
 $$
 \begin{equation} 
 +\Im \int_{\T}|P_N u|^2 |\partial^2_x z|^2 \, dx  \label{yt}
\end{equation}
and notice that the last term vanishes. We thus get  thanks to (\ref{tot}) and (\ref{L5}), 
\begin{eqnarray} 
 I & \lesssim & \||P_N u|^2  z \|_{X^{-1/2+,1}_{]t-\delta ,t+\delta[ }} \|  z \|_{X^{1/2+,1}_{]t-\delta ,t+\delta[ }} 
 + \Bigl\| 2\partial_x( |P_N u|^2)   \partial_x z+
 \partial_x^2( |P_N u|^2)     z \Bigr\|_{X^{-1/2+,0}_{]t-\delta ,t+\delta[ }} \|  z \|_{X^{1/2+,2}_{]t-\delta ,t+\delta[ }} 
  \nonumber \\
  & \lesssim &  M_0^2 \varepsilon(N)  \|  z \|_{X^{1/2+,2}_{]t-\delta ,t+\delta[ }} 
    +\Bigl(M_0^2 N  \varepsilon(N)^{1/2} \|  z \|_{X^{1/2+,2}_{]t-\delta ,t+\delta[ }}^{1/2}+M_0^2 N^2   \varepsilon(N)\Bigr) \|  z \|_{X^{1/2+,2}_{]t-\delta ,t+\delta[ }} \nonumber \\
     & \lesssim &   \varepsilon(N) \|  z \|_{X^{1/2+,2}_{]t-\delta ,t+\delta[ }}^2 +C(N) \; .
\end{eqnarray}
Combining this last estimate with (\ref{estit2})-(\ref{estit9}), (\ref{est54}) follows. 

 We thus infer that
 \begin{eqnarray}
 \| z(t+\delta)\|_{H^2(\T) }^2 &\le & \| z(t)\|_{H^2(\T) }^2 e^{-\gamma \delta}+C(N)\nonumber \\
 & & +
 C_1\, \Bigl( \varepsilon(N)-C_2 (1-e^{-\gamma \delta})\Bigr) \| z(t)\|_{H^2(\T) }^2  \; .
 \end{eqnarray}
 For $ N $ large enough the last term of the right-hand side is clearly negative and is bounded
  from above by
  $$
  - \alpha \, \| z(t)\|_{H^2(\T) }^2\; ,
  $$
  for some small real number $ \alpha>0 $.
 This ensures that, taking $ N>0 $ large enoug,  there exists $ C(N)>0 $ such that
 \begin{equation}
  \|z(t)\|_{H^2(\T) }\le C(N), \; \forall t\ge 0,  \; \label{eestZ'}
 \end{equation}
 and thus on account of (\ref{eq1vZ}), (\ref{estg}), (\ref{estZ'}) and the definition of $ z$, there exists $ K(N)>0 $ such that 
 \begin{equation}
   \|v(t)\|_{H^2(\T) }\le K(N),\; \forall t\ge 0 \; .\label{eestv}
 \end{equation}
 This completes the proof of Proposition \ref{prop2}.
 ההההההההההההההההההההההההההההההההההההההה
   \subsection{Compactness in $ H^2(\T) $}
  To prove the compactness in $ H^2(\T) $, it suffices to show that 
   \begin{equation}
  \|z(t)\|_{H^2(\T) }\le \varepsilon(N),\; \forall t\ge 0 \; .\label{eestzz}
 \end{equation}
 Indeed, this will imply the same estimate on $ v$ and thus on any $ a\in {\mathcal A} $ which will 
 clearly prove the $ H^2(\T) $ compactness of $ {\mathcal A} $.
For proving (\ref{eestzz}), we revisit Lemma \ref{lemestit} with (\ref{boundH2}) at hand.   It is then easy to check that the terms involving $ g_N $ in Lemma 
  \ref{lemestit} (see (\ref{estit2})-(\ref{estit3})) 
    can now be controlled by
 \begin{equation}\label{estcomp1}
 \| g_N \|_{H^2(\T)} \Bigl( \| g_N \|_{H^2(\T)}^2 + \| P_N u\|_{L^\infty(I;H^2(\T))}^2
 + \| z \|_{L^\infty(I;H^2(\T))}^2\Bigr) \lesssim \varepsilon(N) \Bigl(  \| z\|_{X^{1/2+,2}_{I}} +1 \Bigr) \; .
 \end{equation}
  and that (see (\ref{yt}) above)
   \begin{eqnarray} 
  \Bigl| &\Im &\int_t^{t+\delta} e^{-\gamma (t+\delta-s)} \Bigl(|P_N u(s)|^2 z(s),z(s)  \Bigr)_{H^2(\T)}\, ds \Bigr| 
  \lesssim   \| |P_N u|^2 z \|_{L^\infty(]t, t+\delta[;H^1(\T))} \| z\|_{L^\infty(]t, t+\delta[; H^1(\T))} \nonumber  \\
& &   + \Bigl(\|  \partial_x(|P_N u|^2) \partial_x z \|_{L^\infty(]t, t+\delta[; L^2(\T))} +
  \|  \partial_x^2(|P_N u|^2) z \|_{L^\infty(]t, t+\delta[; L^2(\T))}\Bigr) \| z\|_{L^\infty(]t, t+\delta[; H^2(\T))} \nonumber \\
    & & \lesssim   \| P_N u \|_{L^\infty(]t, t+\delta[;H^2(\T))}^2 ( N^{-2} +N^{-1})
   \| z\|_{L^\infty(]t, t+\delta[; H^2(\T))}^2 \nonumber  \\
  & &    \lesssim  N^{-1}   \| z\|_{X^{1/2+,2}_{]t, t+\delta[}} \; . 
    \end{eqnarray}
 To conclude we  need the following estimate that we will prove hereafter.
\begin{equation} \label{estcomp4}
\Bigl\|Q_N\Bigl( |P_N u|^2 P_N u\Bigr) \Bigr\|_{ X^{-1/2+,2}_{]t-1,t+1[}}\lesssim \varepsilon(N),
\; \forall t\ge 0  \; .
\end{equation}
 Proceeding as in the derivation of  (\ref{est54}) but with(\ref{estcomp1})-(\ref{estcomp4}) at hand 
   it is now easy to see that it  actually  holds
 \begin{eqnarray*}
\Bigl|   \Im \Bigl[  \int_t^{t+\delta} e^{-\gamma (t+\delta-s)}\Bigl(Q_N(|v(s)|^2 v(s)) , z(s))\Bigr)_{H^2}\, ds
  \Bigr]\Bigr|
  & \lesssim &  \varepsilon(N) (1+ \| z\|_{X^{1/2+,2}_{]t,t+\delta[}}^2)
\end{eqnarray*}
and thus 
  \begin{eqnarray}
 \| z(t+\delta)\|_{H^2(\T) } &\le & \| z(t)\|_{H^2(\T) } e^{-\gamma \delta}+\varepsilon(N)\nonumber \\
 & & +
 C_1\, \Bigl( \varepsilon(N)-C_2 (1-e^{-\gamma \delta})\Bigr) \| z(t)\|_{H^2(\T) }^2   \; .
 \end{eqnarray}
 which proves  (\ref{eestzz}).
  \subsubsection{\it Proof of Estimate  (\ref{estcomp4}). } 
  Note first that due to the frequency projections it clearly holds 
  \begin{eqnarray*}
  \Bigl\|Q_N\Bigl( |P_N u|^2 P_N u\Bigr) \Bigr\|_{ X^{-1/2+,2}_{]t,t+\delta[}} & \lesssim &
  \Bigl\|Q_N\Bigl( |P_N u|^2 P_N Q_{N/3} u\Bigr) \Bigr\|_{ X^{-1/2+,2}_{]t,t+\delta[}} \\
  &  & +\Bigl\|Q_N\Bigl( (P_N u)^2 \overline{P_N Q_{N/3} u}\Bigr)
   \Bigr\|_{ X^{-1/2+,2}_{]t,t+\delta[}}
  \end{eqnarray*}
  Therefore, on account of (\ref{boundH2}) it is easy to check that we have only to care 
   about
   \begin{equation}\label{gy}
   \Bigl\|Q_N\Bigl( |P_N u|^2  \partial_x^2 P_N Q_{N/3} u\Bigr) \Bigr\|_{ X^{-1/2+,0}_{]t-1,t+1[}}
  +\Bigl\|Q_N\Bigl( (P_N u)^2 \overline{\partial_x^2P_N Q_{N/3} u}\Bigr) \Bigr\|_{ X^{-1/2+,0}_{]t-1,t+1[}}
   \end{equation}
   since  (\ref{estcomp4}) is  obvious for terms that involved less that two derivatives on $Q_{N/3} u $.
   
   To bound (\ref{gy})  we will use  that there exist $ C>0 $ such that for all 
    $ t\ge 0  $,  $ u-g \in X^{1/2+,2}_{]t-1,t+1[}$ (see (\ref{defg}) for the definition of $g $) with 
    \begin{equation} \label{clam}
    \|u-g \|_{X^{1/2+,2}_{]t-1,t+1[}} \le C \; .
    \end{equation}
    Indeed from (\ref{eestZ'}) and (\ref{estZ'}) we know that for $ N_0>0 $ large enough there exists $ \delta>0 $ and 
    $ C(N_0)>0 $ such that 
    $$
    \|Q_{N_0} (v-g)\|_{X^{1/2+,2}_{]t-\delta,t+\delta[}}=  \|z \|_{X^{1/2+,2}_{]t-\delta,t+\delta[}} \le C(N_0) , \quad \forall t\ge 0 , 
    $$
    and thus $ v-g$ is bounded in $  X^{1/2+,2}_{]t-\delta,t+\delta[}$ uniformly in $ t\ge 0 $. Now, proceeding as in the proof of (\ref{boundH2}) 
       we can decompose $ u$  as $ u=v_n+w_n $ with $ \| w_n\|_{L^\infty(]-1,+\infty[; L^2(\T))} \to 0 $ as $ n\to \infty $   and  $\| v_n-g \|_{X^{1/2+,2}_{]t-1,t+1[}}\le C $ for all $ t\ge 0 $. We thus infer that 
     $$
     w_n \rightharpoonup 0  \mbox{ weakly star in } L^\infty(]t-1,t+1[; L^2(\T)) 
     $$
     and thus 
     $$
     v_n -g \rightharpoonup u-g \mbox{ in }  X^{1/2+,2}_{]t-1,t+1[}\; 
     $$
     which proves (\ref{clam}).\\
     Now, on account of Lemma \ref{lemestg}, it clearly holds
   $$
      \Bigl\|Q_N\Bigl( |P_N u|^2  \partial_x^2 P_N Q_{N/3} g\Bigr) \Bigr\|_{ X^{-1/2+,0}_{]t-1,t+1[}}
  +\Bigl\|Q_N\Bigl( (P_N u)^2 \overline{\partial_x^2P_N Q_{N/3} g}\Bigr) \Bigr\|_{ X^{-1/2+,0}_{]t-1,t+1[}}
     $$
     $$
     \lesssim \| Q_{N/3} g \|_{H^2} \| P_N u \|_{L^\infty(\R_+;H^2(\T))}^2 
     \lesssim \varepsilon(N) \; .
    $$
  It thus remains to estimate 
  $$
      \Bigl\|Q_N\Bigl( |P_N u|^2  \partial_x^2 P_N Q_{N/3} (u-g)\Bigr) \Bigr\|_{ X^{-1/2+,0}_{]t-1,t+1[}}
  +\Bigl\|Q_N\Bigl( (P_N u)^2 \overline{\partial_x^2P_N Q_{N/3} (u-g)}\Bigr) \Bigr\|_{ X^{-1/2+,0}_{]t-1,t+1[}}
     $$
 We take extensions $ \theta=P_N \theta  $ of  $  P_N u $ and $ \vartheta= Q_{N/3} P_N \vartheta $ of $ Q_{N/3}P_N(u-g) $ such that 
 $ \|  \theta  \|_{X^{1/2,0}}\le 2 \| P_N u \|_{X^{1/2,0}_{]t-1,t+1[}} $ and 
 $ \|  \vartheta  \|_{X^{1/2,2}}\le 2 \|Q_{N/3} P_N \vartheta  \|_{X^{1/2,2}_{]t-1,t+1[}} $. By duality it suffices   to prove that for
  $ \epsilon>0 $ small enough,
    $$
     \sup_{\|h\|_{X^{1/2-\epsilon,0}}=1} \Bigl[\Bigl| \Bigl(h,Q_N( |\theta|^2 \partial_x^2 \vartheta)
    \Bigr)_{L^2} \Bigr|+ \Bigl| \Bigl(h,Q_N( \theta^2 \partial_x^2 \overline{\vartheta})
    \Bigr)_{L^2} \Bigr|\Bigr]
      \lesssim \varepsilon(N) \| \theta\|_{X^{1/2,0}}^2   \| \vartheta\|_{X^{1/2,2}}\; 
     $$
     and thus   to estimate
     \begin{eqnarray*}
     J: &=&  \int_{\R^3} \sum_{(k_1,k_2,k_3)\in A(N)}
      |\widehat{h}(\tau,k)|
 |\widehat{\theta}(\tau_1,k_1)| |\widehat{\overline{\theta}}(\tau_2,k_2)||k_3|^2
|\widehat{{\vartheta}}(\tau_3,k_3)|\, d\tau_1\, d\tau_2\, d\tau_3 \\
& & +   \int_{\R^3} \sum_{(k_1,k_2,k_3)\in A(N)}
      |\widehat{h}(\tau,k)|
 |\widehat{\theta}(\tau_1,k_1)| |\widehat{\theta}(\tau_2,k_2)||k_3|^2
|\widehat{\overline{\vartheta}}(\tau_3,k_3)|\, d\tau_1\, d\tau_2\, d\tau_3 \\
     \end{eqnarray*}
     where $ \tau=\tau_1+\tau_2+\tau_3 $, $ k=k_1+k_2+k_3 $ and
     $$
     A(N):=\{(k_1,k_2,k_3)\in \Z^3, \, |k_i|\le N \mbox{ for }i\in\{1,2,3\},\, |k_3|> N/3 , \,  N< |k_1+k_2+k_3|\le 3 N \,  \}\; .
     $$
     On $ \R^3\times A(N) $, the resonance relation (\ref{resonant}) clearly yields 
  $$
  \max(|\sigma|,|\sigma_1|,|\sigma_2|,|{\tilde \sigma}_3|)\gtrsim   |(k_1+k_3)(k_2+k_3) |\gtrsim  N^2\; , 
  $$
  where  $ \sigma=\tau+k^2, \, \sigma_1=\tau_1+k_1^2, \, \sigma_2=\tau_2+k_2^2 $
  and $ \tilde{\sigma}_3=\tau_3-k_3^2 $. Moreover, noticing that $ k_1+k_2 \neq 0 $ on $ A(N) $, we infer that 
  $$
  \max(|\sigma|,|\sigma_1|,|{\tilde \sigma}_2|,|\sigma_3|)\gtrsim   |(k_1+k_2)(k_3+k_2) |\gtrsim  N   \; , $$
  where  ${\tilde \sigma_2}=\tau_2-k_2^2 $
  and $ {\sigma}_3=\tau_3+k_3^2 $.
   Therefore proceeding as in the proof of Lemma \ref{damping} we obtain
\begin{eqnarray*}
J  & \lesssim & N^{-1/8+\epsilon}\|h \|_{X^{1/2-\epsilon,0}} \prod_{i=1}^3 \|\theta  \|_{X^{1/2,0}}^2
 \|\vartheta  \|_{X^{1/2-,2}}\\
 & \lesssim & C M_0^2 N^{-1/8+\epsilon}\|h \|_{X^{1/2-\epsilon,0}}\;  ,
\end{eqnarray*}
which completes the proof of  (\ref{estcomp4}).
\section*{Acknowledgements} 
The author is grateful to the Referee for several valuable remarks. 
The author was partially supported by the ANR project "
Etude qualitative des EDP dispersives".

\end{document}